\documentclass[11pt, final]{amsart}
\usepackage[mathscr]{eucal}

\theoremstyle{plain}
\newtheorem{theorem}{Theorem}
\newtheorem{lemma}{Lemma}

\theoremstyle{remark}
\newtheorem{remark}{Remark}

\newcommand{\cA}{\mathscr{A}}

\newcommand{\cF}{\mathscr{F}}
\newcommand{\cH}{\mathscr{H}}
\newcommand{\la}{\langle}
\newcommand{\od}{\odot}
\newcommand{\ot}{\otimes}
\newcommand{\ra}{\rangle}
\newcommand{\tr}{\operatorname{tr}}

\begin{document}
\title[ALGEBRAIC APPROACH TO WIGNER'S THEOREM]
{AN ALGEBRAIC APPROACH TO WIGNER'S UNITARY-ANTIUNITARY THEOREM}
\author{LAJOS MOLN\' AR}
\address{Institute of Mathematics\\
         Lajos Kossuth University\\
         4010 Debrecen, P.O.Box 12, Hungary}
\email{molnarl@math.klte.hu}
\dedicatory{Dedicated to R\'eka Anna}
\thanks{  This research was supported from the following sources:
          (1) Hungarian National Foundation for Scientific Research
             (OTKA), Grant No. T--016846 F--019322,
          (2) A grant from the Ministry of Education, Hungary, Reg.
             No. FKFP 0304/1997}
\maketitle

\vskip 25pt
\SMALL
\begin{center}
\textsc{Abstract}
\end{center}
We present an operator algebraic approach to Wigner's
uni\-ta\-ry-anti\-uni\-ta\-ry theorem using some classical results
from ring theory.
To show how effective this approach is, we prove a generalization
of this celebrated theorem for Hilbert modules over
matrix algebras. We also present a Wigner-type result for maps on prime
$C^*$-algebras.

\vskip 6pt
1991 \textit{Mathematics subject classification (Amer. Math.
Soc.)}: primary 46C05, 46C50; secondary 47D25,
16N60.

\vskip 30pt

\normalsize
\section{Introduction and statement of the results}

Wigner's unitary-antiunitary theorem reads as follows.
Let $H$ be a complex Hilbert space and let $T:H \to H$ be a surjective
map (linearity is not assumed) with the property that
\[
|\la Tx,Ty\ra |=|\la x,y\ra | \qquad (x,y \in H).
\]
Then $T$ is of the form
\[
Tx=\varphi(x)Ux \qquad (x\in H),
\]
where $U:H \to H$ is an either unitary or antiunitary operator (that is,
$U$ is either an inner product preserving linear bijection or a
bijective conjugate-linear map with the property that $\la Ux,Uy\ra
=\la y,x\ra$ for all $x,y\in H$)
and $\varphi :H \to \mathbb C$ is a so-called phase-function which means
that its values are of modulus one.
This celebrated result plays a very important role in
quantum mechanics and in representation theory in physics.

There are several proofs of this theorem in the literature.
See, for example, \cite{LoMe}, \cite{Rat}, \cite{ShAl}, \cite{Uhl} and
the references therein.
The common characteristic of the arguments presented in those papers is
that they manipulate within the Hilbert space which seems to be very
natural, of course.
In this paper we offer a different approach to Wigner's theorem. Namely,
instead of working in $H$, we push the problem to a certain operator
algebra over $H$ and apply some well-known results from ring theory  to
obtain the desired conclusion.
We should remark that in relation to Wigner's theorem, operator
algebras appear also in the papers of Uhlhorn \cite{Uhl} and Wright
\cite{Wri}.
However, in \cite{Uhl} they have nothing to do with the proof of
the unitary-antiunitary theorem. Indeed, Uhlhorn presents an argument
which
can be classified into the first mentioned group of proofs. Moreover,
in \cite{Wri} the author uses Gleason's theorem at a crucial point of
the proof which is a deep result with long proof.
The advantage of our algebraic approach is that in the classical case,
our proof is very clear and short, and it uses only a well-known
theorem of Herstein on Jordan homomorphisms of rings whose proof needs
only few lines of elementary algebraic computation.
It is noteworthy that this result of
Herstein was known before the first complete proofs of Wigner's theorem
appeared.
Furthermore, which
is more important, our approach makes it possible to generalize
Wigner's original theorem for
Hilbert modules, that is, for inner product structures where the inner
product takes its values in an algebra, not necessarily in the complex
field.
Considering the previously mentioned proofs, they are
based on such characteristic properties of the complex field that
one would meet very serious difficulties if one tried to
reach our more general result using those methods (in fact, we are
convinced that such an approach simply cannot be successful).

We now turn to our results. Let $A$ be a $C^*$-algebra.
Let $\cH$ be a left $A$-module with a map $[.,.] : \cH \times \cH \to A$
satisfying
\begin{enumerate}
\item[(i)]      $[f+g,h]=[f,h]+[g,h]$
\item[(ii)]     $[af,g]=a[f,g]$
\item[(iii)]    $[g,f]=[f,g]^*$
\item[(iv)]     $[f,f]\geq 0$ and $[f,f]=0$ if and only if $f=0$
\end{enumerate}
for every $f,g,h\in \cH$ and $a\in A$.
If $\cH$ is complete with respect to the the norm $f \mapsto \|
[f,f]\|^{1/2}$, then we
say that $\cH$ is a Hilbert $A$-module
with generalized inner product $[.,.]$.
This concept is due to Kaplansky \cite{Kap} and in its full generality
to Paschke \cite{Pas}. Nowadays, Hilbert modules over $C^*$-algebras
play a very important role for example in the K-theory of
$C^*$-algebras.

There is another concept of Hilbert modules due to Saworotnow
\cite{Saw}. These are modules over $H^*$-algebras.
$H^*$-algebras are common generalizations of $L^2$-algebras
(convolution algebras) of compact groups and
Hilbert-Schmidt operator algebras on Hilbert spaces.
The only formal difference in the definition is that in the
case of Saworotnow's modules, the generalized inner product takes its
values in the trace-class of the underlying $H^*$-algebra
and the norm with respect to which we require completeness is
$f \mapsto (\tr [f,f])^{1/2}$. Here, $\tr$ denotes the trace-functional
corresponding to $A$ (see \cite{SaFr}). We should note that
Saworotnow originally posed another axiom, namely, a Schwarz-type
inequality \cite[Definition 1]{Saw}. However, as we proved in
\cite[Theorem]{MolPub}, this axiom is redundant. Saworotnow's modules
appear naturally when dealing with multivariate stochastic processes
(see \cite[Section 5]{WiMa}, \cite{Mas}). Moreover, as it turns out
from \cite[Section 3]{Cno}, for example, they have applications in
Clifford analysis
and hence in some parts of mathematical physics. The theory of these
modules is more satisfactory in the sense that many more Hilbert
space-like results have counterparts in Hilbert modules over
$H^*$-algebras than in Hilbert modules over $C^*$-algebras.
Note that it seems to be more common to use right modules instead of
left ones. Of course, this is not a real difference, only a question of
taste.

If $A=M_d(\mathbb C)$ the algebra of all $d\times d$ complex matrices,
then, $A$ being
finite dimensional, the norms on $A$ are all equivalent. Therefore,
the Hilbert modules over the $C^*$-algebra $M_d(\mathbb C)$ are the
same as the Hilbert modules over the $H^*$-algebra $M_d(\mathbb C)$.

Theorem 1 generalizes the original
unitary-anti\-uni\-tary theorem.
As usual, in a $C^*$-algebra $A$,
$|a|$ denotes the absolute value of the element $a$ which is the unique
positive square-root of $a^*a$. If $\cH$ is a Hilbert module, then
the linear
bijection $U:\cH \to \cH$ is called $A$-unitary if $U(af)=aUf$ $(f\in
\cH, a\in A)$ and $[Uf,Uf']=[f,f']$ holds true for every $f,f' \in \cH$.

\begin{theorem}\label{T:wigmod}
Let $\cH$ be a Hilbert module over the matrix algebra $A=M_d(\mathbb C)$
and suppose that there exist vectors $g,h\in \cH$ such that
$[g,h]=I$. Let $T:\cH \to \cH$ be a surjective function
with the property that
\begin{equation}\label{E:wigcond}
|[Tf,Tf']|=|[f,f']| \qquad (f,f' \in \cH).
\end{equation}

If $d>1$, then there exist an $A$-unitary operator
$U:\cH \to \cH$ and a phase-function
$\varphi :\cH \to \mathbb C$ such that
\[
Tf=\varphi(f)Uf \qquad (f\in \cH).
\]

If $d=1$, then there exist an either unitary or antiunitary
operator $U$ on $\cH$ and a phase-function
$\varphi :\cH \to \mathbb C$ such that
\[
Tf=\varphi(f)Uf \qquad (f\in \cH).
\]
\end{theorem}

It seems natural to ask what happens if $A$ is an infinite
dimensional algebra. We have the following result for trivial modules
over prime $C^*$-algebras. If $\cA$ is a $C^*$-algebra, then $\cA$ is a
left module over itself and if we set $[f,g]=fg^*$ $(f,g \in \cA)$, then
$\cA$ becomes a Hilbert module over $\cA$. This is what we mean
when speaking about trivial modules.
A ring $\mathcal R$ is called prime if for any $a,b \in
\mathcal R$, the relation $a\mathcal R b=\{ 0\}$ implies that either
$a=0$ or $b=0$. For example, every algebra of operators which contains
the ideal of all finite rank operators is easily seen to be prime.
Moreover, von Neumann algebras with trivial centre, that is, factors,
are prime $C^*$-algebras.

\begin{theorem}\label{P:wigC*}
Let $\mathscr{A}$ be a prime $C^*$-algebra with unit and let $\phi
:\mathscr{A} \to \mathscr{A}$ be a surjective function such that
\begin{equation}\label{E:wigfact}
        |\phi(A)\phi(B)^*|=|AB^*| \qquad (A,B \in \cA).
\end{equation}
Then there exist a unitary element $U\in \cA$ and a phase-function
$\varphi: \cA \to \mathbb C$ such that $\phi$ is of the form
\[
        \phi(A)=\varphi(A) AU \qquad (A\in \cA).
\]
\end{theorem}

This result is in accordance with Theorem~\ref{T:wigmod}. In fact,
every $A$-linear operator on the trivial module $\cA$ is equal to the
operator of right multiplication by an element of $\cA$. It is easy to
see that if such a map is $A$-unitary, then the corresponding element of
$\cA$ is unitary.

Finally, we give a new proof of the real version of Wigner's theorem.

\begin{theorem}\label{T:wigreal}
Let $H$ be a real Hilbert space and $T:H \to H$ be a surjective function
with the property that
\begin{equation}\label{E:wigi}
|\langle Tx,Ty\rangle |=|\langle x,y\rangle | \qquad (x,y\in H).
\end{equation}
Then there exist a unitary operator $U:H \to H$ and a
function $\varphi :H \to \{ -1, 1\}$ such that $T$ is of the form
\[
Tx=\varphi(x)Ux \qquad (x\in H).
\]
\end{theorem}

The proofs of the results are based on the following theorems
from ring theory:
\begin{itemize}
\item[-] Herstein's homomorphism-antihomomorphism theorem for
Jordan homomorphisms which map onto prime rings.
\item[-] A result of Martindale on elementary operators on prime rings.
\item[-] A theorem of Martindale (or a result of Jacobson and Rickart)
on the extendability of Jordan homomorphisms defined on the symmetric
elements of a ring with involution.
\end{itemize}

\section{Proofs}

As mentioned in the introduction, Saworotnow's modules have
many convenient properties which are familiar in the theory of Hilbert
spaces. First of all, if $\cH$ is a Hilbert module over an
$H^*$-algebra $A$, then $\cH$ is a Hilbert
space with the inner product $\la .,.\ra =\tr [.,.]$. If $M\subset \cH$
is a closed submodule, then for the closed submodule $M^p=\{ f\in \cH \,
:\, [f,g]=0 \, (g\in M)\}$ we obtain $M^p=M^\perp$.
So, we have the orthogonal decomposition $\cH =M \oplus M^p$
\cite[Lemma 3]{Saw}.
A linear operator $T$ on $\cH$ which is bounded with respect
to the Hilbert space norm defined above is called an $A$-linear
operator if $T(af)=aTf$
holds true for every $f\in \cH$ and $a\in A$.
Every $A$-linear operator $T$ is
adjointable, namely, the adjoint $T^*$ of $T$ in the Hilbert space
sense is $A$-linear and we have $[Tf,g]=[f,T^*g]$ $(f,g\in
\cH)$ \cite[Theorem 4]{Saw}.
Consequently, the collection of all $A$-linear operators forms a
$C^*$-subalgebra
of the full operator algebra on the Hilbert space $\cH$.

For the proof of our Theorem~\ref{T:wigmod} we need the following lemma.
In the case of a Hilbert module $\cH$ over an $H^*$-algebra, the
natural equivalent of the Hilbert base is the so-called modular base
\cite{MolCZ}. A family $\{f_\alpha\}_\alpha \subset \cH$ is said to be
modular orthonormal if
\begin{itemize}
\item[(a)]  $[f_\alpha, f_\beta ]=0$ if $\alpha\neq \beta$,
\item[(b)] $[f_\alpha, f_\alpha]$ is a minimal projection in $A$ for
every $\alpha$.
\end{itemize}
A maximal modular orthonormal family of vectors in $\cH$ is called a
modular base.
The common cardinality of modular bases in $\cH$ is called the modular
dimension of $\cH$ (see \cite[Theorem 2]{MolCZ}).

\begin{lemma}\label{L:findim}
Let $\cH$ be a Hilbert $A$-module over the matrix algebra
$A=M_d(\mathbb C)$. If $M\subset \cH$ is a submodule which is
generated by finitely many vectors, then $M$ has finite modular
dimension.
\end{lemma}

\begin{proof}
Observe that since $A$ is finite dimensional, the submodule
generated by finitely many vectors has finite linear dimension.
Therefore, every such submodule is closed.
Let $M$ be generated by the vectors $f_1, \ldots ,f_n$.
Consider the submodule $M_1=Af_1\subset M$. By
orthogonal decomposition we can write $f_2=g_2+h_2$, where $h_2\in
M_1$, $g_2\in M\cap M_1^p$. Clearly, $M_2=Ag_2 \subset
M_1^p$ and we have $f_1, f_2\in M_1+M_2
\subset M$. Next, let $f_3=g_3+h_3$, where $h_3\in M_1+M_2$
and $g_3 \in M\cap (M_1+M_2)^p$.
Let $M_3=Ag_3$. We have $f_1,f_2,f_3 \in M_1+M_2+M_3\subset
M$. Continuing the process we obtain vectors
$g_1,g_2,\ldots ,g_n$ with $[g_i,g_j]=0$ $(i\neq j)$ for which
$f_1,\ldots ,f_n $ is included in the submodule generated
by $g_1,\ldots ,g_n$. Consequently, $M$ is generated by the $g_k$'s.

Let $g\in \cH$ be a nonzero vector. Write $[g,g]=\sum_n \lambda_n^2
e_n$, where the $e_k$'s are pairwise orthogonal minimal projections.
Let $h_k=(1/\lambda_k)e_k g$. Apparently, we have $[h_i,h_j]=0$ $(i\neq
j)$ and $[h_k,h_k]$ is a minimal projection. We assert that
$\sum_k \lambda_k e_k h_k=g$.
This can be verified by taking the generalized inner product of both
sides of this equation with $g$ and then with any vector $f\in \cH$ for
which $[f,g]=0$.
Collecting the $h$'s corresponding to the generating vectors $g_1,\ldots
,g_n$ of $M$, by \cite[Theorem 1]{MolCZ} we obtain a finite modular base
in $M$.
\end{proof}

\begin{remark}
The previous lemma tells us that, under the above assumption on $\cH$,
a submodule of $\cH$ has finite modular dimension if and only if it has
a finite linear dimension.

To emphasize how different the behaviour of Hilbert modules
can be from that of Hilbert spaces, we note
that in general the statement of the previous lemma does not hold
true for Hilbert modules over infinite dimensional $H^*$-algebras.
\end{remark}

In what follows we define operators which are the natural
equivalent of the finite rank operators in the case of Hilbert spaces.
If $f,g \in
\cH$, then let $f\od g$ denote the $A$-linear operator defined by
\[
(f \od g)h=[h,g]f \qquad (h \in \cH).
\]
It is easy to see that for every $A$-linear operator $S$ we
have
\begin{equation}\label{E:prop1}
S(f\od g)=(Sf)\od g ,\qquad (f\od g) S=f\od (S^*g)
\end{equation}
and
\begin{equation}\label{E:prop2}
(f\od g)(f'\od g')=([f',g]f)\od g'=f\od ([g,f']g').
\end{equation}
Define
\[
\cF(\cH) = \{ \sum_{k=1}^n f_k\od g_k  \, : \, f_k, g_k \in \cH\,
(k=1,\ldots ,n),\,\, n\in \mathbb N \}
\]
which is a *-ideal of the $C^*$-algebra of all $A$-linear operators.
We note that if $\cH$ is a Hilbert module over $M_d(\mathbb C)$, then
the range of every element of $\cF(\cH)$ has finite linear dimension,
but
there can be finite rank operators on the Hilbert space $\cH$ which do
not belong to $\cF(\cH)$. In general, if the underlying $H^*$-algebra
is infinite dimensional, then these two classes of operators have
nothing to do with each other.

The following lemma is a spectral theorem for the
self-adjoint elements of $\cF(\cH)$.

\begin{lemma}\label{L:specresol}
Let $\cH$ be a Hilbert module over the matrix algebra $A=M_d(\mathbb
C)$. If $S\in \cF(\cH)$ is a self-adjoint operator,
then $S$ can be written in the form
\[
S=\sum_{k=1}^n \lambda_k f_k \od f_k
\]
where $\{ \lambda_1, \ldots , \lambda_n\}\subset \mathbb R$
and $\{ f_1, \ldots, f_n\} \subset \cH$ is modular orthonormal.
\end{lemma}

\begin{proof}
Let $S\in \cF(\cH)$ be a self-adjoint operator.
Since the range of $S$ has finite linear dimension, $S$ can be written
in the form
\[
S=\sum_k \lambda_k E_k,
\]
where the $\lambda_k$'s are the pairwise different nonzero
eigenvalues of $S$ and the $E_k$'s are the corresponding
spectral projections. Since $S$ is $A$-linear, its eigensubspaces are
submodules. Hence, every spectral projection is $A$-linear with
range included in the range of $S$. Lemma~\ref{L:findim} yields that the
range of $E_k$ has finite modular dimension.
Choose a modular base in the range of every $E_k$. Using
the analog of the Fourier expansion given in
\cite[Theorem 1, (iv)]{MolCZ} we easily conclude that
$S$ can be written in the desired form.
\end{proof}

Now, we are in a position to prove our first theorem. For the proof we
need the concept of Jordan homomorphisms.
A linear map $\phi$ between algebras $\mathcal A$ and $\mathcal B$
is said to be a Jordan homomorphism if it satisfies
\[
\phi(x)^2=\phi(x^2) \qquad (x\in \mathcal A),
\]
or equivalently
\[
\phi(xy+yx)=\phi(x)\phi(y)+\phi(y)\phi(x) \qquad (x,y\in \mathcal A).
\]

\begin{proof}[Proof of Theorem~\ref{T:wigmod}]
We define a linear transformation $\psi$ on the set of all self-adjoint
elements of $\cF(\cH)$ as follows. For any $\{ \lambda_1,
\ldots, \lambda_n \} \subset \mathbb R$ and
$\{ f_1, \ldots, f_n\} \subset \cH$ (we do not require
modular orthonormality) if
\begin{equation}\label{E:welldef}
S=\sum_k \lambda_k f_k \od f_k,
\end{equation}
then let
\[
\psi(S)=\sum_k \lambda_k Tf_k \od Tf_k.
\]
To see that $\psi$ is well-defined, let $\mu_l \in \mathbb R$ and
$g_l\in \cH$ be such that
\[
\sum_k \lambda_k f_k \od f_k=
\sum_l \mu_l g_l \od g_l.
\]
We compute
\begin{gather*}
[(\sum_k \lambda_k Tf_k \od T f_k)Th, Th]=
\sum_k \lambda_k [Th, T f_k][T f_k, Th]=\\
\sum_k \lambda_k [h, f_k][f_k, h]=
[(\sum_k \lambda_k f_k \od  f_k)h, h]=\\
[(\sum_l \mu_l g_l \od  g_l)h, h]=
\sum_l \mu_l [h, g_l][g_l, h]=\\
\sum_l \mu_l [Th, Tg_l][Tg_l, Th]=
[(\sum_l \mu_l Tg_l \od T g_l)Th, Th].
\end{gather*}
Since $T$ is surjective, we obtain that $\psi$ is well-defined.
Due to the fact that in the form \eqref{E:welldef}
of $S$ we have not required anything from the vectors $f_k$, we obtain
readily that $\psi$ is additive and real linear.

We next show that $\psi$ is a Jordan homomorphism. Let
\[
S=\sum_k \lambda_k f_k \od f_k
\]
where $\{\lambda_1, \ldots, \lambda_n\} \subset \mathbb R$ and  $\{
f_1, \ldots, f_n\}\subset \cH$ is modular orthonormal.
If $\{f,g \}\subset \cH$ is modular orthonormal, then
according to \eqref{E:prop2} we have
\[
f\od f \cdot g\od g=0
\]
and
\[
f\od f \cdot f\od f= ([f,f]f)\od f= f\od f
\]
where we have used the equality $[f,f]f=f$ (see \cite[Lemma 1]{MolCZ}).
Therefore, we have $S^2=\sum_k \lambda_k^2 f_k \od f_k$.
Since, as it is easy to see, $\{ Tf_1,\ldots ,Tf_n\}$
is modular orthonormal, we have
$\psi(S)^2=\sum_k \lambda_k^2 Tf_k \od Tf_k$. This results in
\[
\psi(S)^2=\psi(S^2).
\]
Consequently, $\psi$ is a Jordan homomorphism, more precisely, a Jordan
automorphism of the self-adjoint elements of $\cF(\cH)$. Linearizing
the equality above, that is, replacing $S$ by $S+R$ we deduce
\[
\psi(S)\psi(R)+\psi(R)\psi(S)=\psi(SR+RS)
\]
for every self-adjoint $S,R\in \cF(\cH)$. It is now easy to check that
the map $\Psi :\cF(\cH) \to \cF(\cH)$ defined by
\[
\Psi(S+iR)=\psi(S)+i\psi(R)
\]
for every self-adjoint $S,R\in \cF(\cH)$ is
a Jordan *-automorphism of $\cF(\cH)$ which extends $\psi$.

We claim that $\cF(\cH)$ is a prime ring.
Let $S, R$ be $A$-linear
operators such that $S(f \od g)R=0$ holds true for every $f,g\in \cH$.
For an arbitrary $a\in A$ we infer
\[
[Rh,g]a[Sf,h']=[Rh,g][S(af),h']=[S(af\od g)R h, h']=0 \quad (h,h'\in
\cH).
\]
Since $A$ is clearly a prime ring, we obtain that for every $f,g,h,h'$
we have either $[Rh,g]=0$ or $[Sf,h']=0$. This
implies that either $S=0$ or $R=0$ holds true verifying the
primeness of $\cF(\cH)$.

A well-known theorem of Herstein \cite{Her} says that every Jordan
homomorphism onto a prime algebra is either a homomorphism or an
antihomomorphism.
Accordingly, $\Psi$ is either a *-automorphism
or a *-anti\-auto\-mor\-phism of $\cF(\cH)$.
Suppose first that it is a *-automorphism.
Let $g,h\in \cH$ be fixed vectors with the property that $[g,h]=I$.
Define a linear operator $U:\cH \to \cH$ by
\[
Uf=\Psi(f\od g)Th \qquad (f\in \cH).
\]
For any $R\in \cF(\cH)$ we have
\begin{equation}\label{E:eloform}
URf=\Psi(Rf \od g)Th=\Psi(R)\Psi(f\od g)Th=\Psi(R) Uf
\qquad (f\in \cH).
\end{equation}
Using \eqref{E:prop2} and \eqref{E:wigcond} we compute
\begin{gather*}
[Uf,Uf]=[\Psi(g\od f)\Psi(f\od g)Th,Th]=
[\Psi(g\od f \cdot f\od g)Th,Th]=\\
[\Psi( \sqrt{[f,f]}g \od \sqrt{[f,f]}g)Th,Th]=
[T(\sqrt{[f,f]}g) \od T(\sqrt{[f,f]}g)Th,Th]=\\
[Th,T(\sqrt{[f,f]}g)][T(\sqrt{[f,f]}g),Th]=
[h,\sqrt{[f,f]}g][\sqrt{[f,f]}g,h]=\\
[h,g][f,f][g,h]=[f,f].
\end{gather*}
Clearly, $U$ is injective. Moreover,
just as in the case of Hilbert spaces, by polarization we obtain
\begin{equation}\label{E:uni}
[Uf,Uf']=[f,f'] \qquad (f,f'\in \cH).
\end{equation}
To show the surjectivity of $U$ we compute
\begin{gather*}
URg=\Psi(R)Ug=\Psi(R)\Psi(g \od g)Th=\Psi(R)(Tg\od Tg)Th=\\
\Psi(R)([Th,Tg]Tg)= [Th,Tg]\Psi(R)Tg.
\end{gather*}
For an arbitrary $f \in \cH$ we have
$\Psi(R)=f\od Th$ for some $R\in \cF(\cH)$.
Thus the range of $U$ contains the vector
\[
[Th,Tg](f\od Th)(Tg)=[Th,Tg][Tg,Th]f=[h,g][g,h]f=f,
\]
verifying the surjectivity of $U$. Now, by \eqref{E:uni}
it follows that $U$ is $A$-linear and hence an $A$-unitary operator.

From \eqref{E:eloform} we get $\Psi(R)=URU^*$ $(R\in \cF(\cH))$.
Therefore, for every $f\in \cH$ we obtain
\[
Tf\od Tf =\Psi(f \od f)= U(f\od f)U^*=Uf\od Uf.
\]
In view of \eqref{E:wigcond}, this gives us that
\begin{gather*}
[f',f][f,f']=[Tf',Tf][Tf, Tf']=[(Tf \od Tf)Tf',Tf']\\
[(Uf\od Uf)Tf',Tf']=[Tf',Uf][Uf,Tf']= [U^*Tf',f][f,U^*Tf']
\end{gather*}
holds true for every $f,f'\in \cH$. Replacing $f$ by $xf$ $(x\in A)$,
we deduce
\[
[f',f]x^*x[f,f']=[U^*Tf',f]x^*x[f,U^*Tf'].
\]
Since every $x\in A$ is a linear combination of positive elements,
we have
\begin{equation}\label{E:vagy-vagy}
[f',f]x[f,f']=[U^*Tf',f]x[f,U^*Tf'] \qquad (x\in A).
\end{equation}
According to a result of Martindale \cite{Marelem} (see \cite[Lemma
1.3.2]{HerInvo}), if an elementary operator
$x\mapsto \sum_{k=1}^n a_kxb_k$ defined on a prime ring $\mathcal R$ is
identically $0$, then $a_1,\ldots, a_n\in \mathcal R$ are linearly
dependent over the
extended centroid of $\mathcal R$ and the same is true for $b_1, \ldots,
b_n\in \mathcal R$. By the remark after \cite[Proposition 2.5]{Mat}, the
extended
centroid of a prime $C^*$-algebra is just $\mathbb C$ (this remarkable
fact will be used also in the proof of our Proposition~\ref{P:wigC*}).
So, from \eqref{E:vagy-vagy} we get that for every $f,f'\in \cH$
the elements $[f,f']$ and $[f,U^*Tf']$ of $A$ are linearly dependent.
Fix $f'\in \cH$. We know that the linear operators $f\mapsto [f,f']$ and
$f\mapsto [f, U^*Tf']$ are locally linearly dependent. It is
elementary linear algebra to verify that in this case these
operators are (globally) linearly dependent.
Hence, we conclude that for every $f'\in \cH$ there is a
scalar $\varphi(f')$ such that $\varphi(f')f'=U^*Tf'$.
It follows that
\[
Tf=\varphi(f)Uf \qquad (f\in \cH).
\]
Since
\[
\|Tf\|^2=\tr [Tf,Tf]=\tr [f,f]=\tr [Uf,Uf]=\| Uf\|^2 \qquad (f\in \cH),
\]
we obtain that $\varphi$ is a phase-function.

It remains to consider the case when $\Psi$ is *-antiautomorphism.
Just as above, let $g,h\in \cH$ be fixed such that
$[g,h]=I$. Define $U: \cH \to \cH$ by
\begin{equation*}\label{E:deff}
Uf=\Psi(g\od f)Th \qquad (f\in \cH).
\end{equation*}
Clearly, $U$ is a conjugate-linear operator.
Similarly to \eqref{E:eloform}, it is easy to verify that
\begin{equation}\label{E:aha}
U Rf=\Psi(R)^*Uf
\end{equation}
holds true for every $R\in \cF(\cH)$ and $f\in \cH$.
Moreover, just as in the case when $\Psi$ is a *-automorphism,
we obtain
\begin{equation}\label{E:aha2}
[Uf,Uf]=[f,f] \qquad (f\in \cH).
\end{equation}
In particular, $U$ is injective.
By \eqref{E:aha}, using an argument similar to what we have
followed in the first case, one can check
that $U$ is surjective. By the conjugate-linearity of $U$,
\eqref{E:aha2} yields
\[
[Uf,Uf']=[f',f] \qquad (f,f'\in \cH).
\]
Let $f\in \cH$ and define $S=f\od f$. From \eqref{E:aha} we obtain
\begin{gather*}
[\zeta, f][f, \xi]=
[(f\od f) \zeta, \xi]=
[\zeta, S \xi]=
[U(S \xi), U \zeta]=\\
[\Psi(S)^*(U \xi), U\zeta]=
[(Tf \od Tf)(U \xi), U\zeta]=
[U \xi, Tf][Tf, U\zeta]
\end{gather*}
for every $\zeta, \xi \in \cH$. This gives us that
\[
[U^{-1} T\xi,f][f, U^{-1} T\xi]=[T\xi, Tf][Tf, T\xi]=
[\xi, f][f,\xi]
\]
holds true for every $f, \xi \in \cH$.
Fixing $\xi$, just as in the case when $\Psi$ is an automorphism,
we obtain that $[f, \xi]$ and $[f, U^{-1}T \xi]$ are linearly dependent
for every $f$. Therefore, $\xi$ and
$U^{-1}T\xi$ are linearly dependent for every $\xi\in \cH$.
This shows that there exists a phase-function
$\varphi: \cH \to \mathbb C$ such that
\begin{equation}\label{E:tom}
Tf = \varphi(f)Uf \qquad (f \in \cH).
\end{equation}
The proof is now complete in the case when $d=1$.

Suppose that $d>1$. In the antiautomorphic case,
by \eqref{E:tom} we have
\[
|[f,f']|=|[Tf,Tf']|=|[Uf, Uf']|=|[f',f]| \qquad (f,f' \in \cH).
\]
Since there are vectors $g,h \in \cH$ such that $[g,h]=I$, it follows
that $|a|=|a^*|$ holds true for every $a\in A=M_d(\mathbb C)$. As $d>1$,
it is an obvious contradiction, so this case cannot arise.
\end{proof}

\begin{remark}
Observe that if $n=1$, that is, when we have the classical
situation of Hilbert
spaces, our proof is much shorter (see \cite[Theorem
1]{MolJNG} where this case was treated) and uses only
Herstein's
homomorphism-antihomomorphism theorem whose proof needs only few lines
of algebraic computation (see, for example, \cite[Theorem
2.1]{Bre} or \cite[6.3.2 Lemma, 6.3.6 Lemma and 6.3.7 Theorem]{Pal}).

From the proof of Theorem~\ref{T:wigmod} it should be clear why we have
considered modules over matrix algebras. Namely, by the structure
theorem of $H^*$-algebras due to Ambrose \cite{Amb}, the full matrix
algebras are the only unital prime $H^*$-algebras.

One may ask the meaning of the existence of two vectors $g,h \in \cH$
with the property $[g,h]=I$ which has appeared in the formulation of the
theorem
above. We claim that this is equivalent to the requirement that the
modular dimension of $\cH$ is not less than $d$. To see this, let $\{
f_1, \ldots, f_d\}\subset \cH$ be modular orthonormal. Choose
appropriate matrices $a_i \in M_d(\mathbb C)$
such that for the vectors $g_i=a_if_i$ we have
$[g_i,g_i]=a_i[f_i,f_i]a_i^* =e_{ii}$
$(i=1, \ldots , d)$, the standard matrix units. It follows that
$[g_1 +\ldots +g_d,g_1 +\ldots +g_d]=I$.
Now, let the modular dimension of $\cH$ be less than $d$ and choose
a modular base $\{ f_1, \ldots, f_n \} \subset \cH$, where $n<d$. By
\cite[Theorem 1, (v)]{MolCZ} and \cite[Lemma 1]{MolCZ} we have
\[
[g,h] = \sum_{k=1}^n [g,f_k][f_k, h]=
\sum_{k=1}^n [g,f_k][f_k,f_k][f_k, h].
\]
Since $[f_k,f_k]$ is a rank-one projection, we obtain that the rank of
$[g,h]$ is not greater than $n$. This shows that $[g,h]\neq I$ for every
$g,h\in \cH$.
It would be interesting to investigate Wigner's theorem also in these
low-dimensional cases.
\end{remark}

\begin{proof}[Proof of Theorem~\ref{P:wigC*}]
Let $A,B\in \cA$ be arbitrary. Define
\[
\psi(A^*A-B^*B)=\phi(A)^*\phi(A)-\phi(B)^*\phi(B).
\]
To see that $\psi$ is well-defined, let $A',B'\in \cA$ be such that
\[
A^*A-B^*B={A'}^*A'-{B'}^*B'.
\]
For every $S\in \cA$ have
\[
SA^*A S^*-SB^*B S^*=S{A'}^*A' S^* -S{B'}^*B' S^*
\]
and by \eqref{E:wigfact} we deduce
\begin{gather*}
\phi(S)\phi(A)^*\phi(A)\phi(S)^*-
\phi(S)\phi(B)^*\phi(B)\phi(S)^*=\\
\phi(S)\phi(A')^*\phi(A')\phi(S)^*-
\phi(S)\phi(B')^*\phi(B')\phi(S)^* .
\end{gather*}
By the surjectivity of $\phi$ there exists an $S\in \cA$ for which
$\phi(S)=I$. We infer
\[
\phi(A)^*\phi(A)-\phi(B)^*\phi(B)=
\phi(A')^*\phi(A')-\phi(B')^*\phi(B').
\]
Therefore, $\psi$ is a well-defined map on the self-adjoint elements of
$\cA$. Using an argument very similar to what we have just applied,
one can prove that $\psi$ is additive.
We claim that $\psi(I)=I$. Let $\phi(S)=I$.
We have $I=\phi(S)^*\phi(S)=\psi(S^*S)$. Since $\psi$ is
positivity preserving and hence monotone, by the inequality $S^*S\leq \|
S\|^2 I \leq n I$, which holds true for some $n\in \mathbb N$, we infer
that
\begin{equation}\label{E:proj}
I=\psi(S^*S)\leq n \psi(I).
\end{equation}
On the other hand, $\psi(I)$ is a projection. In fact, we have
$\psi(I)=\phi(I)^*\phi(I)$. From \eqref{E:wigfact} it
follows that $\psi(I)^3=\psi(I)^2$. By the spectral mapping theorem this
means that the spectrum of $\psi(I)$ is included in $\{ 0,1\}$. Since
$\psi(I)$ is self-adjoint, we obtain
that $\psi(I)$ is a projection. From \eqref{E:proj} we now get
$\psi(I)=I$. If $U=\phi(I)$, then we have $U^*U=\psi(I)=I$. On
the other hand, by \eqref{E:wigfact} it follows that
$UU^*=\phi(I)\phi(I)^*=II^*=I$. Consequently, $U\in \cA$ is unitary.
From \eqref{E:wigfact} we deduce that
\[
AT^*TA^*=\phi(A)\phi(T)^*\phi(T)\phi(A)^*=\phi(A)\psi(T^*T)\phi(A)^*
\]
holds true for every $A,T\in \cA$. Choosing $A=I$,
we obtain $\psi(T^*T)=U^*(T^*T)U$. Therefore, we get
\[
AT^*TA^*=(\phi(A)U^*)T^*T(\phi(A)U^*)^*.
\]
Since this equation holds true for every $T\in \cA$, it follows that
\begin{equation*}\label{E:centro}
AXA^*=(\phi(A)U^*)X(\phi(A)U^*)^* \qquad (X\in \cA).
\end{equation*}
By the primeness of $\cA$, using \cite[Lemma 1.3.2]{HerInvo}
and the remark after \cite[Proposition 2.5]{Mat} just as in the proof
of Theorem~\ref{T:wigmod}, it follows that
for every $A\in \cA$ the elements $A$ and $\phi(A)U^*$ are linearly
dependent. Consequently,
there exists a scalar valued function $\varphi :\cA \to \mathbb C$
such that
\[
\varphi(A) AU=\phi(A) \qquad (A\in \cA).
\]
This relation yields that
\[
|\varphi(A)|^2 \|A\|^2=\|\phi(A)\|^2=
\| \phi(A)\phi(A)^*\| =\| AA^*\|=\| A\|^2 \quad (A\in \cA)
\]
which implies $|\varphi(A)|=1$ $(0\neq A\in \cA)$.
\end{proof}

\begin{proof}[Proof of Theorem~\ref{T:wigreal}]
If $x,y\in H$, then let $x\ot y$ be the rank-one operator defined by
$(x\ot y)z=\la z,y\ra x$ $(z\in H)$.
By the real version of the spectral theorem, every symmetric (that is,
real self-adjoint) finite rank operator $S$ can be written in the form
\[
S=\sum_{k=1}^n \lambda_k x_k \ot x_k,
\]
where $\lambda_k\in \mathbb R$ and
$x_k\in H$. Similarly to the proof of
Theorem~\ref{T:wigmod}, we define $\psi(S)$ by the formula
\[
\psi(S)=\sum_{k=1}^n \lambda_k T x_k \ot Tx_k.
\]
Repeating the argument in the corresponding part of the proof of
Theorem~\ref{T:wigmod},
we see that $\psi$ is a Jordan automorphism of the
symmetric elements of the ring $F(H)$ of all finite rank operators on
$H$.
In what follows suppose that $\dim H \geq 2$. In fact, if $H$ is
one-dimensional, then the statement of the theorem is trivial.

Consider the unitalized algebra
$F(H)\oplus \mathbb R I$ (of course, we have to add the identity
only in the infinite dimensional case). Defining $\psi(I)=I$, we can
extend $\psi$ to
the set of all symmetric elements of the enlarged algebra in an obvious
way. Now we are in a position to apply two
general algebraic results of Martindale on the extension of Jordan
homomorphisms of the symmetric elements of rings with involution
\cite{Mar}. To be precise, in \cite{Mar} Jordan homomorphism means an
additive map $\phi$ which, besides $\phi(s)^2=\phi(s^2)$, satisfies
$\phi(sts)=\phi(s)\phi(t)\phi(s)$ as well. But if the ring in question
is 2-torsion free (in particular, if it is an algebra), this second
equality follows from the first one (see, for example, the proof of
\cite[6.3.2 Lemma]{Pal}).
The statements \cite[Theorem 1]{Mar} in case $\dim H
\geq 3$ and \cite[Theorem 2]{Mar} when $\dim H =2$ imply that
$\psi$ can be extended uniquely to an associative homomorphism of
$F(H)\oplus \mathbb R I$ into itself. To be honest, since the
results of Martindale
concern rings and hence linearity does not appear,
we could guarantee only the additivity of the extension of $\psi$.
However, the construction in \cite{Mar} clearly shows that in the
case of algebras, linear Jordan homomorphisms have linear
extensions. By the uniqueness of the extension it is apparent
that the extension is *-preserving. Next, observe that our
extension maps $F(H)$ into itself. Thus we have an associative
*-homomorphism $\Psi$ of $F(H)$ into itself which extends $\psi$.
It is easy to see that $\Psi$ is a bijection. Indeed, for arbitrary
nonzero vectors $x,y \in H$ pick a vector $z\in H$ with
$\la x,z\ra , \la y,z\ra \neq 0$.
Plainly, $x\ot y$ is a nonzero scalar multiple of
the operator $x\ot x\cdot z\ot z \cdot y\ot y$. Since our $\psi$ is a
bijection from the set of symmetric elements of $F(H)$ onto itself
and $\Psi$ is an (associative) homomorphism, we obtain that every rank-one
operator is in the range of $\Psi$. This proves the surjectivity of
$\Psi$. The injectivity follows from the simplicity of the algebra
$F(H)$.

The form of *-automorphisms of subalgebras of the full operator algebra
on $H$ containing $F(H)$ is well-known. It follows easily from
\cite[3.2. Corollary]{Che}, for example, that
there is a unitary operator $U$ on $H$ for which
\[
\Psi(A)=UAU^* \qquad (A\in F(H)).
\]
This gives us that
\[
Tx \ot Tx=\Psi(x\ot x)=U (x\ot x)U^*=(Ux)\ot (Ux) \qquad (x\in H).
\]
This implies that
$Tx$ is a scalar multiple of $Ux$, and the scalar must be of
modulus one. The proof is now complete.
\end{proof}

\begin{remark}
Observe that in the case when $n\geq 3$ we could have used a theorem of
Jacobson and Rickart \cite[Theorem 5]{JaRiSymm}. Nevertheless, we
referred
to Martindale's paper since it covers the two-dimensional case as well.
\end{remark}



\begin{thebibliography}{99999}

\bibitem[Amb]{Amb}
W. Ambrose,
'Structure theorems for a special class of Banach algebras',
\emph{Trans. Amer. Math. Soc.} \textbf{57} (1945), 364--386.

\bibitem[Bre]{Bre}
M. Bre\v sar,
'Jordan and Lie homomorphisms of associative rings', preprint.

\bibitem[Che]{Che}
P.R. Chernoff,
'Representations, automorphisms, and derivations of some operator
algebras',
\emph{J. Funct. Anal.} \textbf{12} (1973), 275--289.

\bibitem[Cno]{Cno}
J. Cnops,
\emph{Hurwitz Pairs and Applications of M\"obius Transformations},
Thesis (University of Gent, 1994).

\bibitem[Her1]{Her}
I.N. Herstein,
'Jordan homomorphisms',
\emph{Trans. Amer. Math. Soc.} \textbf{81} (1956), 331--341.

\bibitem[Her2]{HerInvo}
I.N. Herstein,
\emph{Rings with Involution}
(University of Chicago Press, 1976).

\bibitem[JaRi]{JaRiSymm}
N. Jacobson and C.E. Rickart,
'Homomorphisms of Jordan rings of self-adjoint elements',
\emph{Trans. Amer. Math. Soc.} \textbf{72} (1952), 310--322.

\bibitem[Kap]{Kap}
I. Kaplansky,
'Modules over operator algebras',
\emph{Amer. J. Math.} \textbf{75} (1953), 839--853.

\bibitem[LoMe]{LoMe}
J.S. Lomont and P. Mendelson,
'The Wigner unitary-antiunitary theorem',
\emph{Ann. Math.} \textbf{78} (1963), 548--559.

\bibitem[Mar1]{Mar}
W.S. Martindale,
'Jordan homomorphisms of the symmetric elements of a ring with
involution',
\emph{J. Algebra} \textbf{5} (1967), 232--249.

\bibitem[Mar2]{Marelem}
W.S. Martindale,
'Prime rings satisfying a generalized polynomial identity',
\emph{J. Algebra} \textbf{12} (1969), 576--584.

\bibitem[Mas]{Mas}
P. Masani,
\emph{Recent trends in multivariate prediction theory}, in Krishnaiah,
P.R., Multivariate Analysis (Academic Press, 1966).

\bibitem[Mat]{Mat}
M. Mathieu,
'Elementary operators on prime $C^*$-algebras. I',
\emph{Math. Ann.} \textbf{284} (1989), 223--244.

\bibitem[Mol1]{MolPub}
L. Moln\'ar,
'A note on the strong Schwarz inequality in Hilbert $A$-modules',
\emph{Publ. Math. (Debrecen)} \textbf{40} (1992), 323--325.

\bibitem[Mol2]{MolCZ}
L. Moln\'ar,
'Modular bases in a Hilbert $A$-module',
\emph{Czech. Math. J.} \textbf{42} (1992), 649--656.

\bibitem[Mol3]{MolJNG}
L. Moln\'ar,
'Wigner's unitary-antiunitary theorem via Herstein's theorem on
Jordan homomorphisms',
\emph{J. Nat. Geom.} \textbf{10} (1996), 137--148.

\bibitem[Pal]{Pal}
T.W. Palmer,
\emph{Banach Algebras and The General Theory of *-Algebras, Vol. I.},
Encyclopedia Math. Appl. 49 (Cambridge University Press,
1994).

\bibitem[Pas]{Pas}
W.L. Paschke,
'Inner product modules over $B^*$-algebras',
\emph{Trans. Amer. Math. Soc.} \textbf{182} (1973), 443--468.

\bibitem[R\"at]{Rat}
J. R\"atz,
'On Wigner's theorem: remarks, complements, comments, and
corollaries',
\emph{Aequationes Math.} \textbf{52} (1996), 1--9.

\bibitem[Saw]{Saw}
P.P. Saworotnow,
'A generalized Hilbert space',
\emph{Duke Math. J.} \textbf{35} (1968), 191--197.

\bibitem[SaFr]{SaFr}
P.P. Saworotnow and J.C. Friedell,
'Trace-class for an arbitrary $H^*$-algebra',
\emph{Proc. Amer. Math. Soc.} \textbf{26} (1970), 101--104.

\bibitem[ShAl]{ShAl}
C.S. Sharma and D.F. Almeida,
'A direct proof of Wigner's theorem on maps which preserve
transition probabilities between pure states of quantum systems',
\emph{Ann. Phys.} \textbf{197} (1990), 300--309.

\bibitem[Uhl]{Uhl}
U. Uhlhorn,
'Representation of symmetry transformations in quantum mechanics',
\emph{Ark. Fys.} \textbf{23} (1963), 307--340.

\bibitem[WiMa]{WiMa}
N. Wiener and P. Masani,
'The prediction theory of multivariate stochastic processes I.',
\emph{Acta Math.} \textbf{98} (1957), 111-150.

\bibitem[Wri]{Wri}
R. Wright,
'The structure of projection-valued states: A generalization of Wigner's
theorem',
\emph{Int. J. Theor. Phys.} \textbf{16} (1977), 567--573.

\end{thebibliography}
\end{document}